\documentclass[12pt]{article}
\usepackage{amsfonts}
\usepackage{latexsym}
\usepackage{amsmath}
\usepackage{amssymb}
\usepackage{amsthm}
\usepackage{fullpage}
\newtheorem{theorem}{Theorem}[section]
\newtheorem{lemma}[theorem]{Lemma}

\newtheorem{proposition}[theorem]{Proposition}
\newtheorem{corollary}[theorem]{Corollary}
\newtheorem{remark}[theorem]{Remark}
\newtheorem{conjecture}[theorem]{Conjecture}
\newtheorem{oldtheorem}{Theorem A.}
\newtheorem{oldlemma}[oldtheorem]{Lemma A.}
\theoremstyle{definition}

\theoremstyle{remark}
\newtheorem*{note*}{Note}



\begin{document}
\small

\title{\bf Remarks on the conjectured log-Brunn-Minkowski inequality}

\medskip

\author {Christos Saroglou}

\date{}

\maketitle

\begin{abstract}
\footnotesize B\"{o}r\"{o}czky, Lutwak, Yang and Zhang recently 
conjectured a certain strengthening of the Brunn-Minkowski inequality for symmetric convex bodies, 
the so-called log-Brunn-Minkowski inequality. We establish this inequality together with its equality cases 
for pairs of convex bodies that are both unconditional with
respect to some orthonormal basis. Applications of this fact are discussed. Moreover, we prove that the log-Brunn-Minkowski inequality 
is equivalent to the (B)-Theorem for the uniform measure of the cube 
(this has been proven by Cordero-Erasquin, Fradelizi and Maurey for the gaussian measure instead).
\end{abstract}

\section{Introduction}

\hspace*{1.7 em}Let $K,\ L$ be two convex bodies (i.e. compact convex sets with non-empty interior), that contain the origin, in $\mathbb{R}^n$ 
and $\lambda\in(0,1)$. The classical Brunn-Minkowski inequality states that
$$V(\lambda K+(1-\lambda)L)^{1/n}\geq\lambda V(K)^{1/n}+(1-\lambda)V(L)^{1/n} \ ,$$
where $\lambda K+(1-\lambda)L=\{\lambda x+(1-\lambda)y\ | \ x\in K,\ y\in L\}$ is the Minkowski convex combination of $K$ and $L$ (with respect to $\lambda$) and $V(\cdot)=V_n(\cdot)$ denotes the $n$-dimensional volume (i.e. Lebesque measure) functional. The Brunn-Minkowski 
inequality has played an essential role in the development of the Theory of Convex Bodies. We refer to \cite{Ga2}, \cite{MP}, \cite{Sch}  for details and references.

The problem of extending the Brunn-Minkowski theory to the $L^p$-setting has attracted much attention in the previous years (see e.g. \cite{LYZ2}, \cite{LYZ3}). 
For $p>0$, define the $L^p$-convex combination of $K$ and $L$
$$\lambda \cdot K+_p(1-\lambda)\cdot L=\{x\in \mathbb{R}^n\ |\ x\cdot u\leq [\lambda h^p_K(u)+(1-\lambda)h^p_L(u)]^{1/p}, \textnormal{ for all }u\in S^{n-1}\} \ ,
$$
where $S^{n-1}$ denotes the unit sphere and $h_K, \ h_L$ are the support functions of $K,\ L$ respectively. The support function $h_K$ of $K$ is defined as
$$h_K(x)=\max\{x\cdot y \ | \ y\in  K\} \ ,\ x\in\mathbb{R}^n \ .$$
One of its basic properties is that if $K_1$, $K_2$ are convex sets, then $K_1\subseteq K_2$ if and only if $h_{K_1}\leq h_{K_2}$. Note also that if $u$ is a unit vector, $h_K(u)$ is 
the distance from the origin of the supporting hyperplane of $K$, whose outer unit normal is $u$. It is not hard to check that the support function of the convex body
$\lambda \cdot K+_p(1-\lambda)\cdot L$ is the largest support function that is less or equal than $[\lambda h^p_K(u)+(1-\lambda)h^p_L(u)]^{1/p}$.

It was shown in \cite{Lut1} that, for $p\geq 1$
\begin{equation}
V(\lambda\cdot K+_p(1-\lambda)\cdot L)^{\frac{p}{n}}\geq \lambda V(K)^{\frac{p}{n}}+(1-\lambda)V(L)^{\frac{p}{n}}\ .
\label{L^p Brunn-Minkowski}
\end{equation}
The case $p=1$ is indeed the classical Brunn-Minkowski inequality.

Several other extensions of classical inequalities are known to be true. For instance, it was proved in \cite{Lut1} that, for $p\geq 1$, the $L^p$-Minkowski inequality holds:
\begin{equation}
\frac{1}{n}\int_{S^{n-1}}h_K^p(x)h_L^{1-p}(x)dS_L(x)\geq V(K)^{\frac{p}{n}}V(L)^{\frac{n-p}{n}}\ ,
\label{L^p Minkowski}
\end{equation}
where $S_L$ is the surface area measure of $L$ viewed as a measure in $S^{n-1}$, defined by:
$$S_L(\omega)=V\Big( \Bigl\{x\in \textnormal{bd}(L):\exists \ u\in \omega,  \textnormal{ so that }
u\textnormal{ is a normal unit vector for } K \textnormal{ at } x\Big\}\Big)\ . $$
The classical Minkowski inequality (which corresponds to the case $p=1$) and its $L^p$-version are powerful tools for the study of various types of isoperimetric problems (see e.g. \cite{Lut2}, \cite{LYZ1} or again \cite{Sch}).

The case $0<p<1$ seems much harder to deal with. As indicated in \cite{BLYZ}, (\ref{L^p Brunn-Minkowski}) and (\ref{L^p Minkowski}) do not hold for all pairs of 
convex bodies $K$ and $L$. In the same article, the authors deal with the problem of whether (\ref{L^p Brunn-Minkowski}) is true for $0<p<1$ 
and for all symmetric convex bodies $K$, $L$ ($K$ is called symmetric if $K=-K$).

Taking limits as $p\rightarrow 0$ in (\ref{L^p Brunn-Minkowski}), one has:
\begin{equation}
V(\lambda\cdot K+_o(1-\lambda)\cdot L)\geq V(K)^{\lambda}V(L)^{1-\lambda}\ ,\label{Main ineq-Unconditional}
\end{equation}
where
$$\lambda\cdot K+_o(1-\lambda)\cdot L=\{x\in \mathbb{R}^n\ |\ x\cdot u\leq  h^{\lambda}_K(u)h^{1-\lambda}_L(u), \textnormal{ for all }u\in S^{n-1}\} \ ,$$
is the 0-convex combination of $K$ and $L$. The purpose of this note is to make some remarks to the following:
\begin{conjecture}\label{conjecture log B-M}
Inequality (\ref{Main ineq-Unconditional}) is true for all $\lambda\in[0,1]$ and for all symmetric convex bodies $K$, $L$ in $\mathbb{R}^n$.
\end{conjecture}
It was shown in \cite{BLYZ} (although not stated explicitly) that inequalities 
(\ref{L^p Brunn-Minkowski}), (\ref{L^p Minkowski}) are actually equivalent in the class of 
symmetric convex bodies in $\mathbb{R}^n$, for $0<\lambda,p<1$ and that (\ref{Main ineq-Unconditional}) 
would imply (\ref{L^p Brunn-Minkowski}), (\ref{L^p Minkowski}) for all $p>0$ (see also \cite{GHW} for another application of the conjectured log-Brunn-Minkowksi inequality). 
Moreover, the authors proved that Conjecture \ref{conjecture log B-M}
is indeed true in the plane:
\begin{oldtheorem}\label{theorem planar} Let $K$, $L$ be symmetric convex bodies in $\mathbb{R}^2$. Then,
$$
V(\lambda\cdot K+_o(1-\lambda)\cdot L)\geq V(K)^{\lambda}V(L)^{1-\lambda}\ ,
$$
with equality if and only if $K$ and $L$ are dilates or if they are parallelograms with parallel sides.
\end{oldtheorem}

A convex body in $\mathbb{R}^n$ will be called unconditional if it is symmetric with respect to the coordinate 
(with respect to our prefixed orthonormal basis) hyperplanes.
For our purposes, an unconditional convex body $K$ will be called irreducible if it cannot be written as the 
cartesian product of unconditional convex bodies. We are now ready to state one of our main results.
\begin{theorem}\label{Main Theorem-Unconditional}
Let $K$, $L$ be unconditional convex bodies (with respect to the same orthonormal basis) in $\mathbb{R}^n$ and $\lambda\in(0,1)$. Then,
$$
V(\lambda\cdot K+_o(1-\lambda)\cdot L)\geq V(K)^{\lambda}V(L)^{1-\lambda}\ .
$$
Equality holds in the following case: Whenever $K=K_1\times\dots\times K_m$, for some irreducible unconditional convex sets $K_1,\dots ,K_m$, then there exist positive numbers $c_1,\dots,c_m$, such that $$L=c_1K_1\times\dots\times c_mK_m \ .$$
\end{theorem}
Thus, Conjecture \ref{conjecture log B-M} is correct for pairs of convex bodies which are symmetric with respect to the same orthonormal basis. It appears that this inequality was more or less known \cite{BL} \cite{CFM}; however the characterization of the equality cases seems to be a new result and will be critical for Corollary \ref{l^0}.

As a consequence, one can establish the $L^p$-Minkowski and the $L^p$-Brunn-Minkowski inequality, $0<p<1$, mentioned previously, together with the limiting case of the latter as $p\rightarrow  0$ (the so-called log-Minkowski inequality) for unconditional convex bodies. The details of how one can derive the equality cases for these inequalities just from the equality cases in (\ref{Main ineq-Unconditional}) are contained in \cite{BLYZ} and we omit them.
\begin{corollary}\label{log-mink}
Let $K,\ L$ be unconditional convex bodies and $0<p<1$. Then,\\
I. $V(\lambda\cdot K+_p(1-\lambda)\cdot L)^{\frac{p}{n}}\geq \lambda V(K)^{\frac{p}{n}}+(1-\lambda)V(L)^{\frac{p}{n}}$.\\
II. $\frac{1}{n}\int_{S^{n-1}}h_K^p(x)h_L^{1-p}(x)dS_L(x)\geq V(K)^{\frac{p}{n}}V(L)^{\frac{n-p}{n}}$.\\
III. $\int_{S^{n-1}}\log (h_K/h_L)h_LdS_L\geq V(L)\log(V(K)/V(L))$.\\
Equality holds in I and II if and only if $K$ is a dilate of $L$ and in III the equality cases are the same as in Theorem \ref{Main Theorem-Unconditional}.
\end{corollary}
The cone-volume measure of $K$ is defined as $S_o(K,\cdot)=h_K(\cdot)S_ K(\cdot)$.
The logarithmic Minkowski problem asks when a measure on $S^{n-1}$ is the cone-volume measure of a convex body. This problem has recently been solved 
by B\"{o}r\"{o}czky, Lutwak, Yang and Zhang \cite{BLYZ2} in the case of even measures. The authors proved that an even measure on $S^{n-1}$ is the cone-volume measure
of a (symmetric) convex body if and only if it satisfies the so called subspace concentration condition. A measure $\sigma$ on $S^{n-1}$  is said to satisfy the subspace concentration condition
if for any subspace $\xi$ of $\mathbb{R}^n$, the following are true:\\
(i)$\ \sigma(\xi\cap S^{n-1})\leq \frac{1}{n}\sigma(S^{n-1})\dim \xi\ .$\\ 
(ii) If equality holds in (i) for a subspace $\xi'$, then there exists a subspace, complementary to $\xi'$, for which euality holds in (i). \\
It turns out that the subspace concentration condition arises naturally in the study of several other important problems in convex geometry (see e.g. \cite{BLYZ3} \cite{HL}).

What still remains unknown concerning the logarithmic Minkowski problem, is the characterization of the cases for which two convex bodies $K$ and $L$ happen to have the 
same cone-volume measure. 
Note that in the traditional $L^1$-case, the surface area
measure of a convex body determines the body up to translations. Theorem
\ref{Main Theorem-Unconditional} allows us to give a complete answer to the previous open problem in the case when $K$ and $L$ are unconditional. 
One can use Corollary \ref{log-mink} (III) in the same way as in \cite{BLYZ}, Theorem 5.2, to obtain:
\begin{corollary}\label{l^0}
Let $K=K_1\times\dots\times K_m$ and $L$ be unconditional convex bodies, where $K_1,\dots,K_m$, are irreducible convex bodies. If $K$ and $L$ have equal cone-volume 
measures, then $$L=c_1K_1\times\dots\times c_mK_m\ , $$ for some $c_1,\dots,c_m>0$.
\end{corollary}
Our main tool for proving (\ref{Main ineq-Unconditional}) for unconditional convex bodies will be the Pr\`{e}kopa-Leindler inequality, which is a functional form of the Brunn-Minkowski inequality. This is synopsized in the following theorem. We refer to \cite{Ga1} for a good survey on the Brunn-Minkowski and the Pr\`{e}kopa-Leindler inequality.
\begin{oldtheorem}\label{Prekopa-Leindler-Theorem}
Let $\lambda\in (0,1)$ and $f,g,h:\mathbb{R}^n\mapsto \mathbb{R}_+$ be functions with the property that for any $x,y\in \mathbb{R}^n$,
\begin{equation}\label{Prekopa-Leindler-inequality}
h(\lambda x+(1-\lambda)y)\geq f^{\lambda}(x)g^{1-\lambda}(y) \ .
\end{equation}
Then,
$$\int_{\mathbb{R}^n}h(z)dz\geq \Big[\int_{\mathbb{R}^n}f(x)dx\Big]^{\lambda}\Big[\int_{\mathbb{R}^n}g(y)dy\Big]^{1-\lambda}\ .$$
\end{oldtheorem}
In order to derive the equality cases in Theorem \ref{Main Theorem-Unconditional}, 
we will make essential use of 
the following result from \cite{Dub} (see also \cite{Sal}).
\begin{oldlemma}\label{Prekopa-Leindler-equality cases}
If equality holds in (\ref{Prekopa-Leindler-inequality}), then $$f(x)=\frac{\int_{\mathbb{R}^n}fdx}{\int_{\mathbb{R}^n}gdy}q^ng(qx+b) \ ,$$
for some $q>0$, $b\in \mathbb{R}^n$.
\end{oldlemma}
Let $K$ be a symmetric convex body in $\mathbb{R}^n$.
The uniform probability measure of $K$ is defined as follows:
$$\mu_K(A)=\frac{V(K\cap A)}{V(K)} \ .$$
Let $\mu$ be a Borel measure in $\mathbb{R}^n$. We will say that $\mu$ has the (B)-property (resp. weak (B)-property) if for any symmetric convex body $K$ in $\mathbb{R}^n$
and for any positive numbers $t_1,\dots , t_n$ (resp. with $t_1=\dots =t_n$), the function
$$s\mapsto\mu\Big(\textnormal{diag}(t_1^{s},\dots , t_n^{s})K\Big)$$
is log-concave. 
Here, $\textnormal{diag}(c_1,\dots , c_n)$ stands for the diagonal matrix whose diagonal entries are $c_1,\dots ,c_n$. 
It was proven in \cite{CFM} that the gaussian measure posseses this property. This result is known as the (B)-Theorem (conjectured by Banaszczyk, see \cite{La}).
\begin{theorem}\label{theorem reduction}
The following are true:\\
(i) If the uniform probability measure of the $n$-dimensional cube has the (B)-property for every positive integer $n$, then the log-Brunn-Minkowski inequality 
(\ref{Main ineq-Unconditional}) 
holds for all positive integers $n$, for all $\lambda\in(0,1)$, and for all symmetric convex bodies $K$, $L$.\\
(ii) If the log-Brunn-Minkowski inequality holds in dimension $n$, then the uniform probability measure of the $n$-dimensional cube has the (B)-property.
\end{theorem}

\begin{theorem}\label{weak (B)-property}
If the log-Brunn-Minkowski inequality is correct in dimension $n$, then the uniform probability measure of every symmetric convex body in $n$-dimensions has the weak (B)-property.
\end{theorem}

\begin{remark}Combining Theorem A. \ref{theorem planar}, Theorem \ref{theorem reduction} and Theorem \ref{weak (B)-property} one can easily obtain:\\
i) The uniform probability measure of the square (i.e. 2-dimensional cube) has the (B)-property.\\
ii) The uniform probability measure of every symmetric planar convex body has the weak (B)-property. 
This provides an alternative proof for the result proved in \cite{Li}.
\end{remark}
This paper is structured as follows: In Section 2, we establish a Minkowski-type inequality for unconditional convex bodies, 
which will follow from Theorem \ref{Main Theorem-Unconditional} and in our opinion is of some interest. 
Theorem \ref{theorem reduction} and Theorem \ref{weak (B)-property} will be proven in Section 3. 
Inequality (\ref{Main ineq-Unconditional}) for unconditional convex bodies will be proven in Section 4. The equality cases are settled in Section 5.
\section{A Minkowski type inequality}

Before proving our main results, we would like to demonstrate an inequality that would follow as a consequence of Conjecture 1.1.

Define the multi-entry version of the 0-convex combination: If $K_1,\dots ,K_m$ are symmetric convex bodies and $\lambda_1,\dots ,\lambda_m$ are positive numbers summing at 1, then
$$\lambda_1\cdot K_1+_o\dots +_o\lambda_m\cdot K_m:=\{x\in \mathbb{R}^n \ | \ x\cdot u\leq h_{K_1}^{\lambda_1}(u)\dots h_{K_m}^{\lambda_m}(u) \ , \textnormal{ for all }u\in S^{n-1}\} \ .$$
Let us first prove that the log-Brunn-Minkowski inequality would be equivalent to its multi-entry analogue.
\begin{lemma}\label{multi variable LBM-Lemma}
Let $\mathcal{C}$ be a class of symmetric convex bodies, closed under 0-convex combinations. Assume that the log-Brunn-Minkowski inequality (\ref{Main ineq-Unconditional}) holds for all $K,L\in \mathcal{C}$ and for all $\lambda\in[0,1]$. Let $m$ be a positive integer, $K_1,\dots ,K_m\in \mathcal{C}$ and ${\lambda}_1,\dots ,{\lambda}_m$ be positive numbers with ${\lambda}_1+\dots +{\lambda}_m=1$. Then,
\begin{equation}
V(\lambda_1\cdot K_1+_o\dots +_o\lambda_m\cdot K_m)\geq V(K_1)^{\lambda_1}\dots V(K_m)^{\lambda_m} \ .\label{multi variable LBM-ineq}
\end{equation}
\end{lemma}
Proof. We have assumed (\ref{multi variable LBM-ineq}) to be true for $m=2$. We will prove our claim using induction on $m$. Set $\lambda={\lambda}_1+\dots +{\lambda}_{m-1}$ and assume that our assertion is true for the integer $m-1$. Since $$\frac{\lambda_1}{\lambda}+\dots +\frac{\lambda_{m-1}}{\lambda}=1 \ ,$$
the inductive hypothesis states that
\begin{equation}
V\Big(\frac{\lambda_1}{\lambda}\cdot K_1+_o\dots +_o\frac{\lambda_{m-1}}{\lambda}\cdot K_{m-1}\Big)\geq V(K_1)^{\frac{\lambda_1}{\lambda}}\dots V(K_{m-1})^{\frac{\lambda_{m-1}}{\lambda}}\ .
\label{inductive hypothesis}
\end{equation}
It suffices to prove (\ref{multi variable LBM-ineq}) for the integer $m$. We have:
\begin{eqnarray*}
\lambda_1\cdot K_1+_o\dots +_o\lambda_m\cdot K_m\\
&=&\Big\{x\in \mathbb{R}^n\ |\ x\cdot u\leq[h_{K_1}^{\frac{\lambda_1}{\lambda}}(u)\dots h_{K_{m-1}}^{\frac{\lambda_{m-1}}{\lambda}}(u)]^{\lambda}h_{K_m}^{1-\lambda}(u)\Big\}\\
&\supseteq& \{x\in \mathbb{R}^n\ |\ x\cdot u\leq h_K^{\lambda}(u)h_{K_m}^{1-\lambda}(u)\} \ ,
\end{eqnarray*}
where
$$K:=\frac{\lambda_1}{\lambda}\cdot K_1+_o\dots +_o\frac{\lambda_{m-1}}{\lambda}\cdot K_{m-1} \ .$$
Since $K\in\mathcal{C}$, using (\ref{Main ineq-Unconditional}) and (\ref{inductive hypothesis}), we obtain:
\begin{eqnarray*}
V(\lambda_1\cdot K_1+_o\dots +_o\lambda_m\cdot K_m)\\
&\geq &V(K)^{\lambda}V(K_m)^{1-\lambda}\\
&=&V(\frac{\lambda_1}{\lambda}\cdot K_1+_o\dots +_o\frac{\lambda_{m-1}}{\lambda}\cdot K_{m-1})V(K_m)^{\lambda_m}\\
&\geq& \Big[V(K_1)^{\frac{\lambda_1}{\lambda}}\dots [V(K_{m-1})^{\frac{\lambda_{m-1}}{\lambda}}\Big]^{\lambda}V(K_m)^{\lambda_m}\\
&=&V(K_1)^{\lambda_1}\dots V(K_m)^{\lambda_m} \ . \ \ \Box
\end{eqnarray*}
\begin{proposition}\label{multi-minkowski}
Assume that the log-Brunn-Minkowski inequality holds. Let $m$ be a positive integer, $p_1,\dots , p_m$ be positive numbers with $p_1+\dots +p_m=1$ and $K,L_1,\dots ,L_m$ be symmetric convex bodies in $\mathbb{R}^n$. Then,
$$\frac{1}{n}\int_{S^{n-1}}h_{L_1}^{p_1}(x)\dots h_{L_m}^{p_m}(x)dS_K(x)\geq V(L_1)^{\frac{p_1}{n}}\dots V(L_m)^{\frac{p_m}{n}}V(K)^{\frac{n-1}{n}} \ .$$
\end{proposition}
Proof. Since we have assumed that the log-Brunn-Minkowski inequality holds, (\ref{multi variable LBM-ineq}) also holds by the previous lemma. Set $$L:=p_1\cdot L_1+_o\dots +_o p_m\cdot L_m \ .$$
Then, $h_L\leq h_{L_1}^{p_1}\dots h_{L_m}^{p_m}$, thus by the Minkowski inequality,
\begin{eqnarray*}
\frac{1}{n}\int_{S^{n-1}}h_{L_1}^{p_1}(x)\dots h_{L_m}^{p_m}(x)dS_K(x)\\
&\geq &\frac{1}{n}\int_{S^{n-1}}h_{L}(x)dS_K(x)\\
&\geq &V(L)^{\frac{1}{n}}V(K)^{\frac{n-1}{n}}\\
&\geq & V(L_1)^{\frac{p_1}{n}}\dots V(L_m)^{\frac{p_m}{n}}V(K)^{\frac{n-1}{n}} \ . \ \ \Box
\end{eqnarray*}
Notice that for $m=2$ and $L_2=K$, the inequality in Proposition \ref{multi-minkowski} is exactly the $L^{p_1}$-Minkowski inequality, which as already mentioned, is known to be equivalent to the $L^{p_1}$-Brunn-Minkowski inequality. Since its limit case (as $p_1\rightarrow \infty$) is exactly the log-Brunn-Minkowski inequality, Proposition \ref{multi-minkowski} is another equivalent formulation of Conjecture 1.1.
We should also remark here that Proposition \ref{multi-minkowski} holds if we restrict ourselves to a class $\mathcal{C}$ (in the same way as in Lemma \ref{multi variable LBM-Lemma}) of convex bodies, which is closed under 0-convex combinations. The class of unconditional convex bodies is of course such an example. Thus, combining Theorem \ref{Main Theorem-Unconditional} and Proposition \ref{multi-minkowski}, we obtain the following:
\begin{corollary} Let $m$ be a positive integer, $p_1,\dots , p_m$ be positive numbers with $p_1+\dots +p_m=1$, $K$ be a convex body and $L_1,\dots ,L_m$ be unconditional convex bodies in $\mathbb{R}^n$. Then,
$$\frac{1}{n}\int_{S^{n-1}}h_{L_1}^{p_1}(x)\dots h_{L_m}^{p_m}(x)dS_K(x)\geq V(L_1)^{\frac{p_1}{n}}\dots V(L_m)^{\frac{p_m}{n}}V(K)^{\frac{n-1}{n}} \ .$$
\end{corollary}
\section{Reduction to the (B)-property}
Proof of Theorem \ref{weak (B)-property}:

If $K$, 
$L$ are symmetric convex bodies in $\mathbb{R}^n$, and $s, \ t$ are positive numbers, set 
$Q_{\lambda}=K\cap(e^{\lambda t+(1-\lambda)s}L)\ , 
\ \lambda\in (0,1) .$
Let $x\in \lambda\cdot Q_0+_{o}(1-\lambda)\cdot Q_1$, $u\in S^{n-1}$. Then, 
\begin{eqnarray*}
x\cdot u\leq h^{1-\lambda}_{K\cap e^sL}(u)\cdot h^{\lambda}_{K\cap e^tL}(u)\\
&\leq&\min\{h_K(u),h_{e^tL}(u)\}^{\lambda}\cdot\min\{h_K(u),h_{e^sL}(u)\}^{1-\lambda}\\
&\leq&\min\{h_K(u),h_{e^tL}^{\lambda}(u)h_{e^sL}^{1-\lambda}\} \ .
\end{eqnarray*}
This shows that 
$x\in Q_{\lambda}$, thus
\begin{eqnarray*}
V(Q_{\lambda})\geq V(\lambda\cdot Q_0+_o(1-\lambda)\cdot Q_1)\\
&\geq& V(Q_0)^{\lambda}V(Q_1)^{1-\lambda}\\
&=&V(K\cap e^tL)^{\lambda}V(K\cap e^sL)^{1-\lambda} \ . \ \Box
\end{eqnarray*}
The rest of this section is dedicated to the proof of Theorem \ref{theorem reduction}.
\begin{lemma}\label{lemma-argue}Let $K$ be a symmetric convex body and $C_n=[-1/2,1/2]^n$ be the $n$-dimensional cube. 
If $t_1,\dots ,t_n>0$, $s,  t\in \mathbb{R}$, $\lambda\in(0,1)$, the following is true:
$$\Big[\textnormal{diag}\Big(t_1^{\lambda s+(1-\lambda)t},\ \dots \ ,t_n^{\lambda s+(1-\lambda)t}\Big)C_n\Big]\cap K$$
$$\supseteq \lambda\cdot \Big\{ \Big[\textnormal{diag}(t_1^s,\dots , t_n^s)C_n \Big]\cap K\Big\}+_o(1-\lambda)\cdot 
\Big\{ \Big[\textnormal{diag}(t_1^t,\dots , t_n^t)C_n \Big]\cap K\Big\} \ .$$
\end{lemma}
Proof. Set $C_n(\lambda):=\textnormal{diag}\Big(t_1^{\lambda s+(1-\lambda)t},\ \dots \ ,t_n^{\lambda s+(1-\lambda)t}\Big)C_n$ and $Q_{\lambda}:
=C_n(\lambda)\cap K$. Let $x\in \lambda\cdot Q_1+_o(1-\lambda)\cdot Q_0$.
Then, $x\cdot u\leq h_{Q_1}^{\lambda}(u)h_{Q_0}^{1-\lambda}(u)\leq h_K^{\lambda}(u)h_K^{1-\lambda}(u)=h_K(u)$, $u\in S^{n-1}$, hence $x\in K$. Also,
\begin{eqnarray*}
x\cdot e_i\\
&\leq&h_{Q_1}^{\lambda}(e_i)h_{Q_0}^{1-\lambda}(e_i)\\
&\leq&h^{\lambda}_{C_n(1)}(e_i)h^{1-\lambda}_{C_n(0)}(e_i)\\
&=&(t_i^s)^{\lambda}(t_i^t)^{1-\lambda}\\
&=&h_{C_n(\lambda)}(e_i) \ , \ i=1\dots ,n \ .
\end{eqnarray*}
Since $C_n(\lambda)$ is a (coordinate) parallelepiped, it follows that $x$ is contained in $C_n(\lambda)$ as well. 
This shows that $x\in K\cap C_n(\lambda)=Q_{\lambda}$, which gives $\lambda\cdot Q_1+_o(1-\lambda)\cdot Q_0\subseteq Q_{\lambda}.\ \Box$
%
\\
\\
Proof of Theorem \ref{theorem reduction}\\
Assume that the log-Brunn-Minkowski inequality is true in dimension $n$. Let $K$ be a symmetric convex body in 
$\mathbb{R}^n$ and $C_n=[-1/2,1/2]^n$ be the $n$-dimensional cube. For $t_1,\dots ,t_n>0$, $s,  t\in \mathbb{R}$, 
$\lambda\in (0,1)$, one may use the previous lemma to obtain
$$V\Big(\Big[\textnormal{diag}\Big(t_1^{\lambda s+(1-\lambda)t},\ \dots \ ,t_n^{\lambda s+(1-\lambda)t}\Big)C_n\Big]\cap K\Big)$$
$$\geq V\Big(\lambda\cdot \Big\{ \Big[\textnormal{diag}(t_1^s,\dots , t_n^s)C_n \Big]\cap K\Big\}+_o(1-\lambda)\cdot \Big
\{ \Big[\textnormal{diag}(t_1^t,\dots , t_n^t)C_n \Big]\cap K\Big\}\Big)$$
$$\geq V\Big( \Big\{ \Big[\textnormal{diag}(t_1^s,\dots , t_n^s)C_n \Big]\cap K\Big\}\Big)^{\lambda} V\Big(\Big
\{ \Big[\textnormal{diag}(t_1^t,\dots , t_n^t)C_n \Big]\cap K\Big\}\Big)^{1-\lambda} \ .$$
Thus, the function
$$\lambda\mapsto V\Big( \Big[\textnormal{diag}\Big(t_1^{\lambda s+(1-\lambda)t},\ \dots \ ,t_n^{\lambda s+(1-\lambda)t}\Big)C_n\Big]\cap K\Big)$$
$$=(t_1\dots t_n)^{\lambda s+(1-\lambda)t}
V\Big( C_n\cap \Big[\textnormal{diag}\Big(s_1^{\lambda s+(1-\lambda)t},\ \dots \ ,s_n^{\lambda s+(1-\lambda)t}\Big)K\Big]\Big) $$ 
is log-concave, where $s_i=t_i^{-1}$. Since $s_1,\dots , s_n$ are arbitrary, it is clear that the uniform probability measure of $C_n$ has the (B)-property.
It remains to prove that (ii) implies (i).

Let $m\geq n$ be an integer, $r_1,\dots ,r_m$, $s_1,\dots , s_m$ be positive numbers and $v_1,\dots , v_m$
be unit vectors. 
Define the set
$$R_{\lambda}=\bigcap_{i=1}^m\Big\{x\in \mathbb{R}^n\ \Big|\ |x\cdot v_i|\leq r_i^{\lambda}s_i^{1-\lambda}\Big\} \ .$$
First observe that, by an approximation argument, if one proved the inequality
$$V(R_{\lambda})\geq V(R_0)^{1-\lambda}V(R_1)^{\lambda} \ ,$$
for all $m$, $r_i$, $s_i$, $i=1,\dots ,m$, the log-Brunn-Minkowski inequality would follow. Indeed, one can choose the sequence of sets 
$\{v_i:i=1,\dots ,m\}\subseteq S^{n-1}$, $m\in \mathbb{N}$, so that the sequence of the uniform discrete probability measures supported on these sets converges weakly
to the uniform probability measure of $S^{n-1}$, as $m\rightarrow \infty$. Then, for any $\lambda\in [0,1]$, $R_{\lambda}$ converges to $\lambda\cdot K+_o(1-\lambda)\cdot L$, with
respect to the Haussdorff metric, as $m\rightarrow \infty$. 

In other words, we need to show that if $$F(\lambda):=\int_{\mathbb{R}^n}\prod_{i=1}^m
\mathbf{1}_{[-r_i^{\lambda}s_i^{1-\lambda},r_i^{\lambda}s_i^{1-\lambda}]}(x\cdot v_i)dx \ ,$$
then
\begin{equation}
F(\lambda)\geq F(0)^{1-\lambda}F(1)^{\lambda} \ .\label{F-ineq}
\end{equation}
Let $\varepsilon>0$. Define the function
$$G_{\varepsilon}(\lambda)=\int_{x\in\mathbb{R}^n}
\int_{u\in \mathbb{R}^m}\prod_{i=1}^m\mathbf{1}_{[-r_i^{\lambda}s_i^{1-\lambda},r_i^{\lambda}s_i^{1-\lambda}]}(u_i+x\cdot v_i)\mathbf{1}_{[-\varepsilon, \varepsilon]}(u_i)du dx\ .$$
Then, by a linear change of variables, we have:
$$G_{\varepsilon}(\lambda)=
\int_{x\in\mathbb{R}^n}\int_{w\in \mathbb{R}^m}
\prod_{i=1}^m\mathbf{1}_{[-r_i^{\lambda}s_i^{1-\lambda},r_i^{\lambda}s_i^{1-\lambda}]}(w_i)\mathbf{1}_{[-\varepsilon, \varepsilon]}(w_i-x\cdot v_i)dw dx\ .$$
Consider the symmetric $(n+m)$-dimensional convex bodies
$$K=\{(x,w)\ |\ x\in \mathbb{R}^n,w\in\mathbb{R}^m, \ |w_i-x\cdot v_i|\leq \varepsilon ,\ i=1,\dots m,\ \|x\|_{\infty}\leq b/2\} \ ,$$
$$C=C_{n+m}=\{(x,w)\ |\ x\in \mathbb{R}^n,w\in\mathbb{R}^m,\ |x_i|\leq 1/2, \ |w_j|\leq 1/2,\ i=1,\dots ,n, \ j=1,\dots,m\} \ ,$$
the $(n+m)$-dimensional cube of volume 1, where $b$ is some positive number. Define also the diagonal map
$$T_{\lambda}=\textnormal{diag}(b/2,\dots,b/2,2r_1^{\lambda}s_1^{1-\lambda},\dots,2r_m^{\lambda}s_m^{1-\lambda}) \ .$$Now, choose $b$ so large that the cylinder 
$\{x\in \mathbb{R}^{n+m} :  \|x\|_{\infty}\leq b/2\}$ contains the $(m+n)$-dimensional convex body defined by the inequalities $|w_i|\leq r_i^{\lambda}s_i^{1-\lambda}$, $|w_i-x\cdot v_i|\leq \varepsilon$,
$i=1,\dots ,m$. It is then clear that
\begin{eqnarray*}
G_{\varepsilon}(\lambda)=V(K\cap T_{\lambda}C)\\
&=&\det T_{\lambda}V((T^{-1}_{\lambda}K)\cap C)\\
&=& 2^{m-n}b^m\Big[\prod_{i=1}^mr_i^{\lambda}s_i^{1-\lambda}\Big]V(T^{-1}_{\lambda}K\cap C)\ .
\end{eqnarray*}
Note also that $$T^{-1}_{\lambda}=\textnormal{diag}(1,\dots ,1,(s_1/r_1)^{\lambda},\dots ,(s_m/r_m)^{\lambda})\cdot \textnormal{diag}((b/2)^{-1},\dots ,(b/2)^{-1},(2s_1)^{-1},
\dots ,(2s_m)^{-1}) \ .$$Thus, if we set
$$K':=\textnormal{diag}((b/2)^{-1},\dots ,(b/2)^{-1},(2s_1)^{-1},
\dots ,(2s_m)^{-1})K \ ,$$
then
$$G_{\varepsilon}(\lambda)=\mu_C\Big(\textnormal{diag}(1,\dots ,1,(s_1/r_1)^{\lambda},\dots ,(s_m/r_m)^{\lambda})K'\Big) \ ,$$
which by assumption (ii) is a log-concave function of $\lambda$.

For $u\in \mathbb{R}^m$, introduce
$$F_u(\lambda)=\int_{\mathbb{R}^n}\prod_{i=1}^m
\mathbf{1}_{[-r_i^{\lambda}s_i^{1-\lambda},r_i^{\lambda}s_i^{1-\lambda}]}(x\cdot v_i+u_i)dx \ .$$
Then, $$G_{\varepsilon}(\lambda)=\int_{u\in[-\varepsilon ,\varepsilon]^m}F_u(\lambda)du \ .$$
Note that $F(\lambda)=F_0(\lambda)$, so by continuity if (\ref{F-ineq}) does not hold, there exists an $\varepsilon >0$, such that for all $u\in[-\varepsilon ,\varepsilon]^m$,
$$F_u(\lambda)<F_u^{\lambda}(1)F_u^{1-\lambda}(0) \ .$$
Using H\"{o}lder's inequality, we obtain:
\begin{eqnarray*}
G_{\varepsilon}(\lambda)=\int_{u\in[-\varepsilon ,\varepsilon]^m}F_u(\lambda)du\\
&<&\int_{u\in[-\varepsilon ,\varepsilon]^m}F_u^{\lambda}(1)F_u^{1-\lambda}(0)du\\
&\leq&\Big[\int_{u\in[-\varepsilon ,\varepsilon]^m}F_u(1)du\Big]^{\lambda}\Big[\int_{u\in[-\varepsilon ,\varepsilon]^m}F_u(0)du\Big]^{1-\lambda}\\
&=&G_{\varepsilon}^{\lambda}(1)G_{\varepsilon}^{1-\lambda}(0) \ .
\end{eqnarray*}
This is a contradiction, since $G_{\varepsilon}(\lambda)$ is a log-concave function of $\lambda$, so (ii) implies (i). $\Box$
\section{The use of the Pr\`{e}kopa-Leindler inequality}

Following  \cite{CFM}, if $K$, $L$ are unconditional convex sets and $\lambda\in[0,1]$, we define
$$K^{\lambda}\cdot L^{1-\lambda}:=\{(\pm |x_1|^{\lambda}|y_1|^{1-\lambda},\dots ,\pm |x_n|^{\lambda}|y_n|^{1-\lambda})\ | \ (x_1,\dots ,x_n)\in K, \ (y_1,\dots ,y_n)\in L\} \ .$$
The key to the proof of (\ref{Main ineq-Unconditional}) for unconditional convex bodies will be
the following simple observation.
\begin{lemma}\label{Main observation}
If $K$, $L$ are unconditional convex bodies and $\lambda\in[0,1]$, then
$$\lambda\cdot K+_o(1-\lambda)\cdot L\supseteq K^{\lambda}\cdot L^{1-\lambda} \ .$$
\end{lemma}
Proof. Since $\lambda\cdot K+_o(1-\lambda)\cdot L$ and $K^{\lambda}\cdot L^{1-\lambda}$ 
are unconditional, it suffices to prove that
$$\lambda\cdot (K\cap\mathbb{R}_+^n)+_o(1-\lambda)\cdot 
(L\cap\mathbb{R}_+^n)\supseteq (K^{\lambda}\cdot L^{1-\lambda})\cap\mathbb{R}_+^n\ .$$
Note that since $\lambda\cdot K+_o(1-\lambda)\cdot L$ is convex and unconditional, one can write
$$\lambda\cdot (K\cap\mathbb{R}^n_+)+_o(1-\lambda)\cdot (L\cap \mathbb{R}^n_+)=
\{x\in \mathbb{R}^n_+\ | \ x\cdot u\leq h_K^{\lambda}(u)h_L^{1-\lambda}(u),
\textnormal{ for all }u\in S^{n-1}\cap \mathbb{R}^n_+\} \ .$$
Let $x\in K\cap\mathbb{R}_+^n$, $y\in L\cap\mathbb{R}_+^n$, $u\in S^{n-1}\cap \mathbb{R}_+^n$. 
Set $$X=((x_1u_1)^{\lambda},\dots ,(x_nu_n)^{\lambda}), \ Y=((y_1u_1)^{1-\lambda},
\dots ,(y_nu_n)^{1-\lambda}) \in \mathbb{R}^n\ .$$Then,
\begin{eqnarray*}
(x_1^{\lambda}y_1^{1-\lambda},\dots ,x_n^{\lambda}y_n^{1-\lambda})\cdot u= X\cdot Y \\
&\leq& \| X \|_{\frac{1}{\lambda}}\|Y\|_{\frac{1}{1-\lambda}}\\
&=& (x\cdot u)^{\lambda}(y\cdot u)^{1-\lambda}\\
&\leq &h_K^{\lambda}(u)h_L^{1-\lambda}(u) \ ,
\end{eqnarray*}
which shows that $(x_1^{\lambda}y_1^{1-\lambda},\dots ,x_n^{\lambda}y_n^{1-\lambda})\in (K^{\lambda}\cdot L^{1-\lambda})\cap\mathbb{R}_+^n$, proving our claim. $\Box$\\
\\
Inequality (\ref{Main ineq-Unconditional}) for unconditional convex bodies, follows immediately from the previous lemma and Proposition 10 from \cite{CFM} 
which states that the volume of $K^{\lambda}\cdot L^{1-\lambda}$ is a log-concave function of $\lambda$, where $K,\ L$ are unconditional convex bodies. 
However, since we want to investigate equality cases, we will need to repeat the proof of the previously mentioned fact.
\begin{proposition}\label{L=TK}
Let $K$, $L$ be unconditional convex bodies in $\mathbb{R}^n_+$ and $\lambda\in(0,1)$. Then,
\begin{equation}
V(K^{\lambda}\cdot L^{1-\lambda})\geq V(K)^{\lambda}V(L)^{1-\lambda} \ . \label{Maurey combination}
\end{equation}
If equality holds in the last inequality, then there exists a positive definite diagonal matrix $T$ with $L=TK$.
\end{proposition}
Proof. For $x\in \mathbb{R}^n$, set $f(x)=\mathbf{1}_{K\cap\mathbb{R}^n_+}(x)$, 
$g(x)=\mathbf{1}_{L\cap\mathbb{R}^n_+}(x)$, $h(x)=\mathbf{1}_{( K^{\lambda}\cdot L^{1-\lambda})\cap\mathbb{R}^n_+}(x)$. 
Set also, $\overline{f}(x)=f(e^{x_1},\dots , e^{x_n})e^{x_1+\dots +x_n}$, 
$\overline{g}(x)=g(e^{x_1},\dots , e^{x_n})e^{x_1+\dots +x_n}$, $\overline{h}(x)=h(e^{x_1},\dots , e^{x_n})e^{x_1+\dots +x_n}$. It is clear by the definitions that $$\overline{h}(\lambda x+(1-\lambda)y)\geq \overline{f}(x)^{\lambda}\overline{g}(y)^{1-\lambda}\ , \ x,y\in \mathbb{R}^n \ .$$Thus, by the Pr\`{e}kopa-Leindler inequality, we obtain:
\begin{equation}
\int_{\mathbb{R}^n}\overline{h}(z)dz\geq\Big[\int_{\mathbb{R}^n}\overline{f}(x)dx\Big]^{\lambda}
\Big[\int_{\mathbb{R}^n}\overline{g}(y)dy\Big]^{1-\lambda} \ .\label{In proposition PL-ineq}
\end{equation}
Next, in the same spirit as in \cite{FM}, Proposition 1, use the change of variables $z_i\mapsto e^{z_i}$, $i=1,\dots ,n$, to obtain:
\begin{eqnarray*}
V(K^{\lambda}\cdot L^{1-\lambda})=2^n\int_{\mathbb{R}^n_+}h(z)dz\\
&=&2^n\int_{\mathbb{R}^n}\overline{h}(z)dz\\
&\geq&\Big[2^n\int_{\mathbb{R}^n}\overline{f}(x)dx\Big]^{\lambda}
\Big[2^n\int_{\mathbb{R}^n}\overline{g}(y)dy\Big]^{1-\lambda}\\
&=&\Big[2^n\int_{\mathbb{R}_+^n}f(x)dx\Big]^{\lambda}
\Big[2^n\int_{\mathbb{R}_+^n}g(y)dy\Big]^{1-\lambda}\\
&=&V(K)^{\lambda}V(L)^{1-\lambda}
\end{eqnarray*}
and (\ref{Maurey combination}) is proved. Now, equality holds in (\ref{Maurey combination}) if and only if 
equality holds in (\ref{In proposition PL-ineq}). According to Lemma A. \ref{Prekopa-Leindler-equality cases}, 
if equality holds in
(\ref{In proposition PL-ineq}), there exist $c>0, \ q\in \mathbb{R}$ and $b\in \mathbb{R}^n$ such that
\begin{equation}
\overline{f}(x)=c\overline{g}(qx+b) \ .\label{f=g(Tx)}
\end{equation}
Notice that the set $$\{x\in \mathbb{R}^n\ | \ (e^{x_1},\dots ,e^{x_n})\in K \ , \ (e^{qx_1+b_1},\dots ,e^{qx_n+b_n})\in L\}$$
has non-empty interior. This fact together with the last inequality imply that inside some open subset of $\mathbb{R}^n$,
$$e^{x_1+\dots +x_n}=c'e^{q(x_1+\dots +x_n)} \ ,$$
for some other constant $c'>0$. This clearly shows that $q=1$, so (\ref{f=g(Tx)}) becomes
$$\mathbf{1}_{K\cap\mathbb{R}^n_+}(e^{x_1},\dots ,e^{x_n})=\mathbf{1}_{L\cap\mathbb{R}^n_+}(e^{b_1}e^{x_1},
\dots ,e^{b_n}e^{x_n}) \ , \ x\in\mathbb{R}^n \ ,$$
which shows that $L=TK$, where $T=\textnormal{diag}(e^{b_1},\dots , e^{b_n})$. $\Box$
\section{Equality cases}
Let $\lambda\in (0,1)$ and $K$, $L$ be convex bodies, such that 
$$V(\lambda\cdot K+_o(1-\lambda)\cdot L)=V(K)^{\lambda}V(L)^{1-\lambda} \ .$$ 
Then, $V(K^{\lambda}\cdot L^{1-\lambda})=V(K)^{\lambda}V(L)^{1-\lambda}$ by Lemma 4.1 and Proposition 4.2, and 
hence $L=TK$ for some positive definite diagonal map $T$. One can easily check that 
$K^{\lambda}\cdot (TK)^{1-\lambda}\supseteq T^{1-\lambda}K$, so it follows by Lemma 4.1 that
$$T^{1-\lambda}K\subseteq \lambda\cdot K+_o(1-\lambda)\cdot (TK)\ .$$
On the other hand, we have $$V(T^{1-\lambda}K)=V(K)^{\lambda}V(TK)^{1-\lambda}=V(\lambda\cdot K+_o(1-\lambda)\cdot (TK)) \ ,$$
thus 
\begin{equation*}\label{main-ineq-section 5}
T^{1-\lambda}K=\lambda\cdot K+_o(1-\lambda)\cdot (TK)\ .
\end{equation*}

\begin{lemma}\label{main-section 5}
Let $T$ be a positive definite diagonal matrix and  $K$ be an unconditional convex body 
such that $$\lambda\cdot K+_o(1-\lambda)\cdot (TK)\subseteq T^{1-\lambda}K \ , $$
for some $\lambda\in(0,1)$. Then, the restriction of $T$ on every irreducible component of $K$ is a multiple of the 
identity.
\end{lemma}
Proof. Suppose that $K=K_1\times\dots\times K_m$, where $K_i$ is an unconditional compact convex set,
$i=1,\dots ,m$. Trivially, 
$$\lambda\cdot (K_1\times \{0\})+_o(1-\lambda)\cdot (TK_1\times\{0\})\subseteq 
\lambda\cdot K+_o(1-\lambda)\cdot (TK) \subseteq T^{1-\lambda}K\ .$$
Thus, $$\lambda\cdot K_1+_o(1-\lambda)\cdot (TK_1)\subseteq 
T^{1-\lambda}K_1 \ .$$
Thus, $K_1$ satisfies the assumption of our lemma, therefore it suffices to assume that $K=K_1$, i.e. $K$ is itself irreducible. We need to prove that $T$ is
a multiple of the identity. Suppose that it is not. Since $K$ is irreducible, it is true that there exists a unit vector 
$v=(v_1,\dots,v_n)$ in the support of the surface area measure of $K$ and a smooth boundary point $x$ which has $v$ as its exterior unit 
normal vector, 
so that none of them lie in a proper coordinate subspace. In other
words, $x_i,v_i\neq 0$, $i=1,\dots,n$. By the fact that $K$ is unconditional, we may actually take all the $x_i$'s and $v_i$'s 
to be positive. Note that $T^{1-\lambda}x$ and $Tx$ are smooth boundary points of $T^{1-\lambda}K$ and $TK$ respectively. Also, 
since $T$ is 
not a multiple of the identity, the exterior unit normal $T^{-1}v$ to $TK$ at $Tx$ is not parallel to $v$. In particular,
$h_{TK}(v)>(Tx)\cdot v$ and $h_K(u)>x\cdot u$, if $u\in S^{n-1}\setminus \{v\}$. It follows by the H\"{o}lder inequality that if 
$u\in (\mathbb{R}^n_+\cap  S^{n-1})\setminus \{v\}$, then
$$u\cdot T^{1-\lambda}x\leq (u\cdot x)^{\lambda}(u\cdot Tx)^{1-\lambda}< h_K^{\lambda}(u)h_{TK}^{1-\lambda}(u)\ .$$
Therefore, $T^{1-\lambda}x$ is an interior point of $\lambda\cdot K+_o(1-\lambda)\cdot (TK)$, which contradicts our assumption. We arrived 
at a contradiction because we assumed that $T$ is not a multiple of the identity. $\Box$
\\

To end the proof of Theorem \ref{Main Theorem-Unconditional}, it remains to prove that when the cases described as "equality cases'' in Theorem \ref{Main Theorem-Unconditional} occur, then equality holds indeed in (\ref{Main ineq-Unconditional}). In particular, since 
it is easily verified that equality holds in (\ref{Main ineq-Unconditional}) for any 0-convex combination of dilates of the same convex body, it suffices to show the 
following: If $\lambda\in [0,1]$, $K_i$ is a convex body in $\mathbb{R}^{k_i}$, $c_i>0$, $i=1,\dots ,m$, 
then $$V\big(\lambda\cdot (K_1\times\dots\times K_m)+_o(1-\lambda)\cdot(c_1K_1\times\dots\times c_mK_m)\big) $$
$$\leq V\big([\lambda\cdot K_1+_o (1-\lambda)\cdot (c_1K_1)]\times\dots\times[\lambda\cdot K_m+_o (1-\lambda)\cdot (c_mK_m)]\big) \ .$$
To see that the last inequality is true, take $x\in \lambda\cdot (K_1\times\dots\times K_m)+_o(1-\lambda)\cdot(c_1K_1\times\dots\times c_mK_m)$. If $u\in \{0_{\mathbb{R}^{k_1+\dots + k_{i-1}}}\}\times \mathbb{R}^{k_i}\times \{0_{\mathbb{R}^{k_{i+1}+\dots +k_m}}\}$, then
$x\cdot u\leq h^{\lambda}_{k_i}(u)(c_ih_{k_i}(u))^{1-\lambda}$, which shows that the projection of $x$ on the subspace $\{0_{\mathbb{R}^{k_1+\dots + k_{i-1}}}\}\times \mathbb{R}^{k_i}\times \{0_{\mathbb{R}^{k_{i+1}+\dots +k_m}}\}$ is contained in $\{0_{\mathbb{R}^{k_1+\dots + k_{i-1}}}\}\times [\lambda\cdot K_i+_o (1-\lambda)\cdot (c_iK_i)]\times \{0_{\mathbb{R}^{k_{i+1}+\dots +k_m}}\}$, $i=1,\dots ,m$. This gives:
$$\lambda\cdot (K_1\times\dots\times K_m)+_o(1-\lambda)\cdot(c_1K_1\times\dots\times c_mK_m)$$
$$\subseteq[\lambda\cdot K_1+_o (1-\lambda)\cdot (c_1K_1)]\times\dots\times[\lambda\cdot K_m+_o (1-\lambda)\cdot (c_mK_m)] \ ,$$
proving our last assertion. $\Box$
\\
\\
\textbf{Acknowledgement.} I would like to thank the referee(s) for 
finding a mistake in Theorem \ref{theorem reduction} and for simplifing significantly the proof of Lemma 5.1.

\bigskip

\noindent \textsc{Ch.\ Saroglou}: Department of Mathematics,
Texas A$\&$M University, 77840 College Station, TX, USA.

\smallskip

\noindent \textit{E-mail:} \texttt{saroglou@math.tamu.edu \ \&\ christos.saroglou@gmail.com}

\end{document}